\documentclass[letterpaper, 10 pt, conference]{ieeeconf}
\IEEEoverridecommandlockouts                              %
\usepackage{algorithm}
\usepackage[noend]{algpseudocode}
\usepackage[table]{xcolor}
\usepackage{graphicx}
\usepackage{epsfig}
 \usepackage{amsmath}
\usepackage{amssymb}
\usepackage{amsfonts}
\usepackage{amsthm}
 \usepackage{epstopdf}
\usepackage{algorithm}
\usepackage{algpseudocode}
\usepackage{cite}

\usepackage{wrapfig}
\usepackage{amsmath}

\usepackage{caption}
\renewcommand{\thesubsection}{\Alph{subsection}}

\newcommand{\mzzz}{\color{black}}
\newcommand{\sy}{\color{black}}
\newcommand{\mz}{\color{black}}

\title{\LARGE \bf
Learning the tuned liquid damper dynamics by means of a robust EKF}

\author{Alberta Longhini, Michele Perbellini, Stefano Gottardi, Shenglun Yi, Hao Liu and Mattia Zorzi \thanks{S. Yi  is with the School of Automation, Beijing Institute of Technology, Beijing 100081, China; A. Longhini, M. Perbellini, S. Gottardi and M. Zorzi is with the Department of Information Engineering, University of Padova, Via Gradenigo 6/B, 35131 Padova, Italy; Liu Hao is with College of Architecture and Civil Engineering, Beijing University of Technology Emails: {\tt\small alberta.longhini@studenti.unipd.it}, {\tt\small michele.perbellini@studenti.unipd.it}, {\tt\small stefano.gottardi@studenti.unipd.it}, {\tt\small 3120185460@bit.edu.cn}, {\tt\small zorzimat@dei.unipd.it}, {\tt \small liuhao\_vv@email.s.bjut.edu.cn}}
}

\begin{document}

\maketitle
\thispagestyle{empty}
\pagestyle{empty}

\begin{abstract}  The tuned liquid dampers (TLD) technology is a feasible and cost-effective seismic design. In order to improve its efficiency it is fundamental to find accurate models describing their dynamic. A TLD system can be modeled through the Housner model and its parameters can be estimated by solving a nonlinear state estimation problem. We propose a robust {\mzzz extended Kalman filter} which alleviates the model discretization and the fact that the noise process is not known. We test the effectiveness of the proposed approach by using some experimental data corresponding to two classical seismic waves, namely the El Centro wave and the Hachinohe wave.
\end{abstract}

\section{Introduction}
It is estimated that about 50,0000 detectable earthquakes occur every year, and some of them are extremely destructive, like the Sichuan earthquake in 2008 or the Aceh earthquake in 2012. Earthquakes can cause great loss of life and buildings, most of which come from damage to houses. The tuned liquid dampers (TLD) technology is a feasible and cost-effective seismic design. A TLD is a passive mechanical damper placed on the top of a building for suppressing the structural vibrations. More precisely,  it is a water tank whose damping effect depends on the sloshing of shallow liquid. Clearly, this technology is also easy to install in existing structures and applicable for temporary use. In order to design an efficient TLD we need to have an accurate model describing it. Housner proposed a nonlinear model for TLD using the analysis of the dynamic behavior for the response of elevated water tanks to earthquake ground motion \cite{housner1963dynamic}. Then,  numerical models have been introduced to solve the equations of the liquid motion \cite{limin1991semi, wakahara1992suppression}. Since then, TLD is widely used in the flexible and weakly-damped structures such as high-rise buildings, towers and suspension bridges \cite{kamgar2020modified, an2019vibration, pandit2020seismic}.

In this paper we consider the problem to learn the parameters of the Housner model for a TLD system from the acceleration and the force applied to the bottom of the tank. Such a problem can be {\sy cast} as a nonlinear state estimation problem, \cite{ding2015simultaneous}. A natural way to tackle this problem is to use the Extended Kalman (EKF) filter, \cite{reif1999stochastic,wang2019second,barrau2016invariant}. In the practical TLD application, however, the implementation of the EKF requires the design of the noise variances which is {\sy lacking relative} prior knowledge. Moreover,  the experimental data is collected with a certain sampling time, so the original continuous model needs to be discretized in order to be used in the EKF. Accordingly, there is a mismatch between the discretized model and the actual  one.

Risk sensitive filtering aims to address the model uncertainty \cite{H_INF_HASSIBI_SAYED_KAILATH_1999, RISK_WHITTLE_1980,levy2016contraction}. More precisely, compared with the traditional way, the standard quadratic loss function is replaced by an exponential quadratic loss function so that large errors are severely penalized.  This severity is tuned by the so called risk sensitive parameter, \cite{boel2002robustness, hansen2005robust, yoon2004robust,OPTIMAL_SPEYER_FAN_BANAVAR_1992}. A refined version of such a paradigm is the robust Kalman filter proposed in \cite{ROBUST_STATE_SPACE_LEVY_NIKOUKHAH_2013} as well as its extensions \cite{zorzi2018robust,abadeh2018wasserstein, zorzi2019distributed, emanuele2020robust, zenere2018coupling,9303925,nguyen2019bridging}. Here, the uncertainty is expressed incrementally at each time step $t$. More precisely, the state estimator is designed according to the least favorable model at time $t$ belonging to the ambiguity set which is a ball about the nominal model in the Kullback-Leibler (KL) topology. Although there are some robust versions of EKF, see \cite{kallapur2009discrete, sahraoui2017robust,kim2020robust}, nobody generalized this paradigm for EKF.

The contribution of this paper is to propose a robust extended Kalman filter, where the uncertainty is expressed incrementally, to learn the TLD dynamics. We adopt the real-time hybrid simulation technology, while two classical seismic waves, namely the El Centro wave and the Hachinohe wave are employed as our input earthquake signals. The aim of this work is to show how the robust algorithm is able to produce a better parameter estimation than the standard one. It turns out that the estimated parameters with the robust approach {\sy tend} not to deviate from the reference parameter (i.e. the actual parameters) keeping the relative error of the estimated parameter very small.

The outline of the paper is as follows. The proposed robust extended Kalman filter is introduced in Section \ref{sec_2}. Section \ref{sec: TDL_learning} introduces the problem to learn the parameters of the Housner model as well as its nonlinear state estimation formulation. Then, the experimental results based on the seismic waves are presented in Section \ref{sec_4}. Finally, Section \ref{sec_5} regards the conclusions and future research directions.

 \section{Robust extended Kalman filter} \label{sec_2}
 We consider the following discrete time state space model:
 \begin{equation}\label{nomi_mod}
     \begin{cases} x_{t+1}= f(x_t,u_t) + B v_t\\
 y_t = h(x_t,u_t) +  D v_t\\
 \end{cases}
 \end{equation}
 where $x_t\in\mathbb{R}^n$ is the state process, $u_t \in \mathbb{R}^q$ is a known deterministic input, and $v_t \in \mathbb{R}^m$ is white Gaussian noise (WGN) with $\mathbb{E}[v_tv_s^{T}] = I_m\delta_{t-s} $, where $\delta_s$ is the Kronecker delta function.
 Our aim is to estimate, in a recursive way, the state $x_{t+1}$ from the observation process $y_t \in \mathbb{R}^p$. If the functions $f(\cdot)$ and $g(\cdot)$ are linear the problem has a well known solution which is given by the \textit{Kalman Filter}. On the contrary, if the two functions, or one of these, are nonlinear we can use the EKF, which aims to find an approximate solution of the problem:
 \begin{equation}\label{min_NL}
     \text{argmin}\, \mathbb{E}\left[||x_{t+1}-g_t(y_t)||^2|Y_{t-1}\right]
 \end{equation}
 where $Y_{t-1} = \left \{y_s, 0 \le s \le t-1 \right\}$.
 This filter adopts at each step a linearization of the nominal model (\ref{nomi_mod}) around the previous estimate. It is well known that the evolution of the estimates of the problem are:
 \begin{align}\label{EKF_eq}
 \hat{x}_{t|t} &={\hat{x}_{t}} + L_t(y_t - h(\hat{x}_t,u_t)) \\
 \hat{x}_{t+1} & = f(\hat{x}_{t|t}, u_t) \notag
 \end{align}
 where $\hat x_{t|t}$ and $\hat x_{t}$ represent the estimate of $x_t$ given $Y_t$ and $Y_{t-1}$, respectively; $L_t$ is the filtering gain.
  However, in the real scenarios the nominal model does not match the actual one. This is due to two main facts. First, the model parameters are affected by  uncertainty, like the real structure of $f(\cdot)$ and the variances of the noises. The second one is that physical models have a natural description using continuous time models. Thus, to use them we have to perform a discretization step, and so what we obtain is an approximation. Therefore we propose a robust approach to compute the gain in (\ref{EKF_eq}) that takes into account the presence of uncertainties. Consider model (\ref{nomi_mod}) and, as in EKF, let us consider {\mz that the state equation is linearized with respect to $\hat x_{t|t}$  while the measurement equation is linearized with respect to $\hat x_t$:
  \begin{equation}\label{linearized}
     \begin{cases} x_{t+1}= A_tx_t - A_t \hat x_{t|t} + f(\hat{x}_{t|t}, u_t) + B v_t\\
    y_t = C_tx_t -C_t \hat x_{t} +  h(\hat{x}_t,u_t)+ D v_t\\
 \end{cases}
 \end{equation}
where $A_t = \partial f(x,u_t)/\partial x|_{x=\hat x_{t|t}}$,  $C_t = \partial h(x,u_t)/\partial x |_{x=\hat x_{t}}$. Let $z_t = \left[x_{t+1}^{T} \ y_t^{T}\right]^{T}$. Then, (\ref{linearized}) is characterized by the nominal transition probability of $z_t$ given $x_t$:
 \begin{align}\label{phi}
     &\phi_t(z_t|x_t) \\
     &\sim \mathcal{N}\left(\begin{bmatrix}A_t x_t - A_t \hat x_{t|t} +f(\hat{x}_{t|t}, u_t) \\ C_t x_t-C_t \hat x_{t} + h(\hat{x}_t,u_t)  \end{bmatrix} , \begin{bmatrix}B \\ D \end{bmatrix}\begin{bmatrix}B^{T} & D^{T}\end{bmatrix}\right) .
 \end{align}
} We assume that the noise $v_t$ affects all the components of the dynamic observations in (\ref{nomi_mod}) and (\ref{linearized}), therefore the covariance matrix:
 \begin{equation}
    K_{z|x} = \begin{bmatrix} B \\ D \end{bmatrix}\begin{bmatrix} B^{T} & D^{T} \end{bmatrix} \notag
 \end{equation}
 is positive definite.
To derive the robust filter we adopt the minimax approach on the linearized model (\ref{linearized}) proposed in the recent works \cite{STATETAU_2017, ROBUST_STATE_SPACE_LEVY_NIKOUKHAH_2013, zorzi2017convergence}. Let $\tilde{\phi}_t(z_t|x_t)$ denote the actual transition probability density of $z_t$ given $x_t$. It is worth noting that $\tilde \phi_t$ is not necessarily Gaussian distributed. Then, we measure the mismatch between the probability densities $\phi_t$ and $\tilde \phi_t$ using the KL divergence:
 \begin{equation}
 \begin{aligned}
     \mathcal{D}(\tilde \phi_t, \phi_t) = \iint \tilde \phi_t (z_t|x_t)& {\mzzz p_t (x_t|Y_{t-1})} \\
    & \ln\left(\frac{ \phi_t(z_t|x_t)}{\tilde{\phi}_{t}(z_t|x_t)}\right)d z_t d x_t.
     \end{aligned}
 \end{equation}
 Then, we assume that $\tilde{\phi}_t$ belongs to $\mathcal{B}_t$, i.e. a ball about the nominal density, with
 \begin{equation}
     \mathcal{B}_t = \left\{\ \tilde{\phi}(z_t|x_t) \ \text{s.t.} \ \mathcal{D}(\tilde{\phi}_t,{\mzzz \phi_t}) \le c_t   \right\}
 \end{equation}
 where $c_t >0$ is the \textit{tolerance} specified at each time step. The latter  measures the mismodeling budget between the nominal and the actual model at time $t$ in terms of KL divergence.
The main idea is to compute the gain $L_t$ in Equation (\ref{EKF_eq}) in a robust way by considering the minimax optimization problem:
\begin{equation}\label{minmax_pb}
    \hat{x}_{t+1} = \text{arg}\min_{g_t \in G_t} \max_{\tilde{\phi}_t \in \mathcal B_t}\mathbb{\tilde E}\left[||x_{t+1}-g_t(y_t)||^2|Y_{t-1}\right]
\end{equation}
where
\begin{equation}\label{def_obj}
\begin{split}
    \mathbb{\tilde E}&\left[||x_{t+1}-g_t(y_t)||^2|Y_{t-1}\right] = \\ &\iint ||x_{t+1} - g_t(y_t)||^2\tilde{\phi}(z_t|x_t){\mzzz p_t(x_t|Y_{t-1})} \ dx_t dz_t
\end{split}
\end{equation}
denotes the mean square error of the estimator $\hat{x}_{t+1} = g_t(y_t)$ of $x_{t+1}$ evaluated with respect to the density $\tilde{\phi}_t$; $G_t$ is the set of estimators for which (\ref{def_obj}) is bounded for any $\tilde \phi_t \in \mathcal B_t$. Therefore, here the gain is computed using as reference the \textit{least favorable} description, and not the nominal one as in Problem (\ref{min_NL}). Under the assumption that the \textit{a priori} probability density of $x_t$ conditioned on the observations $Y_{t-1}$ is:
\begin{equation}
  {\mzzz  p_t(x_t|Y_{t-1})} \sim \mathcal{N}(\hat{x}_t , V_t),
\end{equation}
then the solution to Problem {\mzzz (9)} is a Kalman like filter, see \cite[Theorem 1]{ROBUST_STATE_SPACE_LEVY_NIKOUKHAH_2013} and \cite{robustleastsquaresestimation,STATETAU_2017,OPTIMALITY_ZORZI}. More precisely, the least favorable solution to (\ref{minmax_pb}) is Gaussian with the same mean
of $\phi_t$ defined in (\ref{phi}), while its variance changes. Since the least favorable (linearised) model is Gaussian, the estimator solution to (\ref{minmax_pb}) is the Bayes estimator and thus the estimator corresponding to the original model is (\ref{EKF_eq}) where the filtering gain  now is
\begin{equation}
    L_t = V_tC_t^{T}(C_tV_tC_t^{T}+DD^{T})^{-1}
\end{equation}
with:
\begin{align}
    V_{t+1} &= (P_{t+1}^{-1} - \theta_tI)^{-1}\nonumber\\
    P_{t+1}&=A_t V_t A_t^T\nonumber\\
    & \hspace{0.2cm}-A_tV_tC_t^T (C_tV_tC_t^T+DD^T)^{-1} C_tV_tA_t^T+BB^T\nonumber
\end{align}
where $\theta_t>0$ is the unique solution to $\gamma(P_{t+1},\theta_t) = c_t$ and
\begin{equation}
    \gamma(P,\theta) = \frac{1}{2}\left\{\log\det\left(I -\theta P) + \text{tr}\left[(I-\theta P)^{-1} - I\right]\right)\right\} \notag.
\end{equation}
It is worth noting that the filtering gain is different from the standard one because the variance of $\tilde \phi_t$ is different from the one of $\phi_t$.
Thus, the idea is to use this \textit{robust} approach to find the gain in the linearized version of model in (\ref{nomi_mod}). This leads to the {\mz robust extended Kalman filter (REKF)}, reported in Algorithm 1. {\sy It is worth noting that the full procedure is similar to the EKF, except for the presence of $V_t$, which requires (at step 8) the computation of $\theta_t$. The latter computation is done in a numerical way, since a closed form solution does not exist.}
\begin{algorithm}
    \caption{Robust Extended Kalman filter at time t}\label{euclid}
    \hspace*{\algorithmicindent} \textbf{Input} $\hat{x}_{t}$, $V_{t}$, $c_t$, $y_t$, $u_t$\\
    \hspace*{\algorithmicindent} \textbf{Output} $\hat x_{t|t}$, $\hat{x}_{t+1}$, $V_{t+1}$
    \begin{algorithmic}[1]

    \State $C_{t} = \frac{\partial h(x,u_t)}{\partial x}\Bigr|_{x = \hat{x}_{t}}$
    \State $L_t = V_{t}C_{t}^{T}(C_{t}V_{t}C_{t}^{T} + DD^T)^{-1}$
    \State $\hat{x}_{t|t} = \hat{x}_{t} + L_t[y_t - h(\hat{x}_{t},u_t)]$
    \State $A_t = \frac{\partial f(x,u_t)}{\partial x}\Bigr|_{x = \hat{x}_{t|t}}$
    \State $\hat{x}_{t+1} = f(\hat{x}_{t|t}, u_t)$
    \State $P_{t+1} = A_tV_{t}A_t^{T}$
    \State \hspace{1cm}$ - A_tV_{t}C_t^{T}(C_{t}V_{t}C_{t}^{T} + D D^T)^{-1}C_tV_{t}A_t^{T} + B B^T$
    \State $\textit{find $\theta_t$ s.t.}\ \gamma(P_{t+1}, \theta_t) = c_t$
    \State $V_{t+1} = (P_{t+1}^{-1} - \theta_t I) ^{-1}$
    \end{algorithmic}
    \label{alg: robust EKF}
\end{algorithm}

 \section{Learning the TLD dynamic} \label{sec: TDL_learning}
TLD is a very important technology nowadays: it
 attenuates the effect of vibrations in the building at a low price.
 Learning the model dynamics of TLD's is then crucial
 for future design needs. A TLD is identified
 by the Housner model \cite{housner1963dynamic},
 which describes the relative displacement of water (here called $d$) in the following way:
 \begin{equation}
 \ddot{d} + 2\xi \omega \dot{d} + \omega^2d = - u
 \label{eq:inpu_outo}
 \end{equation} where $u$ is the acceleration on the bottom of the water tank; $\xi$ is the damping ratio; $\beta$ is the ratio between the mass of water that can oscillate horizontally against a restraining spring and the total mass of water in the tank $m_t$; $\omega$ is the frequency of oscillation of water.
The goal is to estimate the parameters of the model,
 through the measurement of the {\mz reactive} force applied to the bottom of the water tank:
 \begin{equation}
 F = -(1-\beta)m_tu+m_t \beta \omega d+ m_t \beta \omega \xi \dot{d}.
 \label{eq: misur}
 \end{equation}
  The total mass of water $m_t$ is known. The parameter $\xi$ is considered as a known constant
 because it is possible to {\sy tune} it as we wish using baffles into the tank.

Thus, the parameters to estimate are $\beta$ and $\omega$. The latter can be treated as states components. Indeed, Housner model in (\ref{eq:inpu_outo})-(\ref{eq: misur}) can be rewritten as the state space model
 \begin{align}
 \dot{x} &=\begin{bmatrix}\dot x_1\\\dot x_2\\ \dot x_3\\ \dot x_4\end{bmatrix}= \begin{bmatrix}-u - 2 x_5 x_4 x_1-x_4^2 x_2\\x_1\\0\\0\end{bmatrix} + v,  \label{eq_state_space1}\\
  y &= -(1-x_3)m_t u+m_t x_3 x_4 x_2+ m_t x_3 x_4 \xi x_1 + w
 \label{eq_state_space2}
 \end{align}
 where:
 \begin{align*}
 	 x_1 = \dot{d}, \quad x_2 = d, \quad
 	x_3 = \beta, \quad
 x_4 = \omega.
 \end{align*} Accordingly, the problem becomes to estimate the state of the nonlinear state model in (\ref{eq_state_space1})-(\ref{eq_state_space2}).
Here, we assume the process noise $v$ and the measurement noise $w$ are Gaussian white noises with variances $Q$ and $R$, respectively. The process variance $Q$ is a diagonal matrix such that:
 \begin{equation}
 Q = diag(\sigma_1, \sigma_2, \sigma_3, \sigma_4).
 \end{equation}
Note that we choose the first two components of the diagonal (corresponding to $x_1$ and $x_2$) with high magnitude to ensure variability over time, while the other two components (relative to $x_3$ and $x_4$) are selected with small value so that the estimated parameters change a little over time. In plain words $\sigma_3$ and $\sigma_4$ tune the a priori information about the changing rate of $x_3$ and $x_4$, respectively, in the stochastic hypermodel for  $\beta$ and $\omega$ in (\ref{eq_state_space1}), \cite{niedzwiecki2000identification}. 	
 Then, the initial state $x_0$ is modeled as a Gaussian random vector with mean $\hat x_0$ and variance $V_0$:
 \begin{equation}
  V_0 = diag(\lambda_1, \lambda_2, \lambda_3, \lambda_4).
 \end{equation}

Finally, it is necessary to discretize the state equations (\ref{eq_state_space1})-(\ref{eq_state_space2}) because the data is collected with a certain sampling time. Here, the fourth-order Runge-Kutta method is used. Clearly this procedure induces an error that deviates the obtained model from the actual one. Furthermore, the implementation of the EKF requires the design of the noises $v$, $w$ that nobody knows how they are really distributed in practice. In other words, there is some model uncertainty between the nominal model and actual model. In order to overcome this problem we propose to use the robust version of the EKF presented in the previous section.


 \section{Experimental results} \label{sec_4}
\begin{figure}[h]
          \begin{center}
            \includegraphics[scale = 0.4]{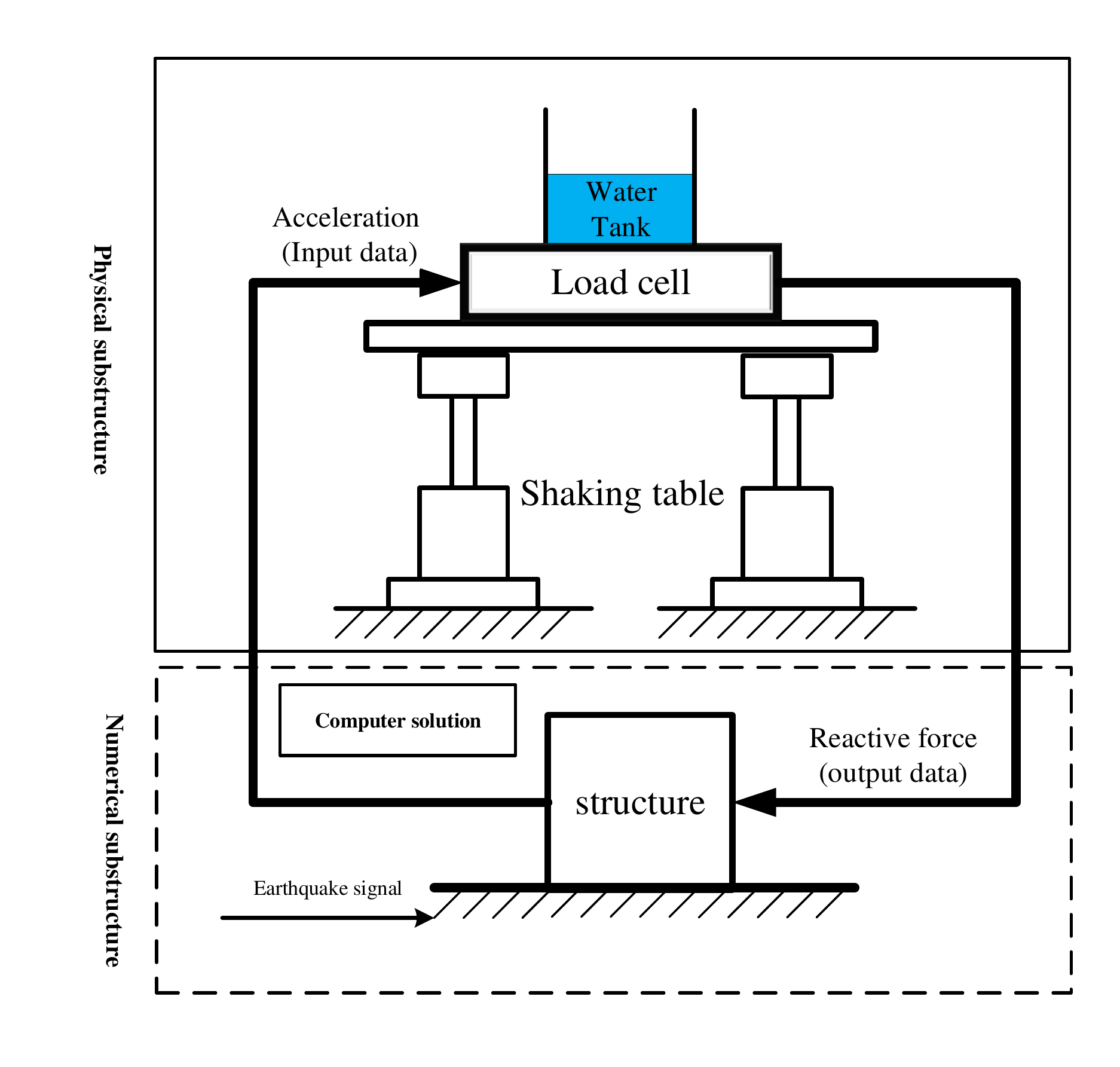}
          \end{center}
          \caption{\footnotesize{Conceptual view of real-time hybrid simulation technology.}}
          \label{fig_tld}
    \end{figure}

\begin{figure}[h]
          \begin{center}
            \includegraphics[scale = 0.35]{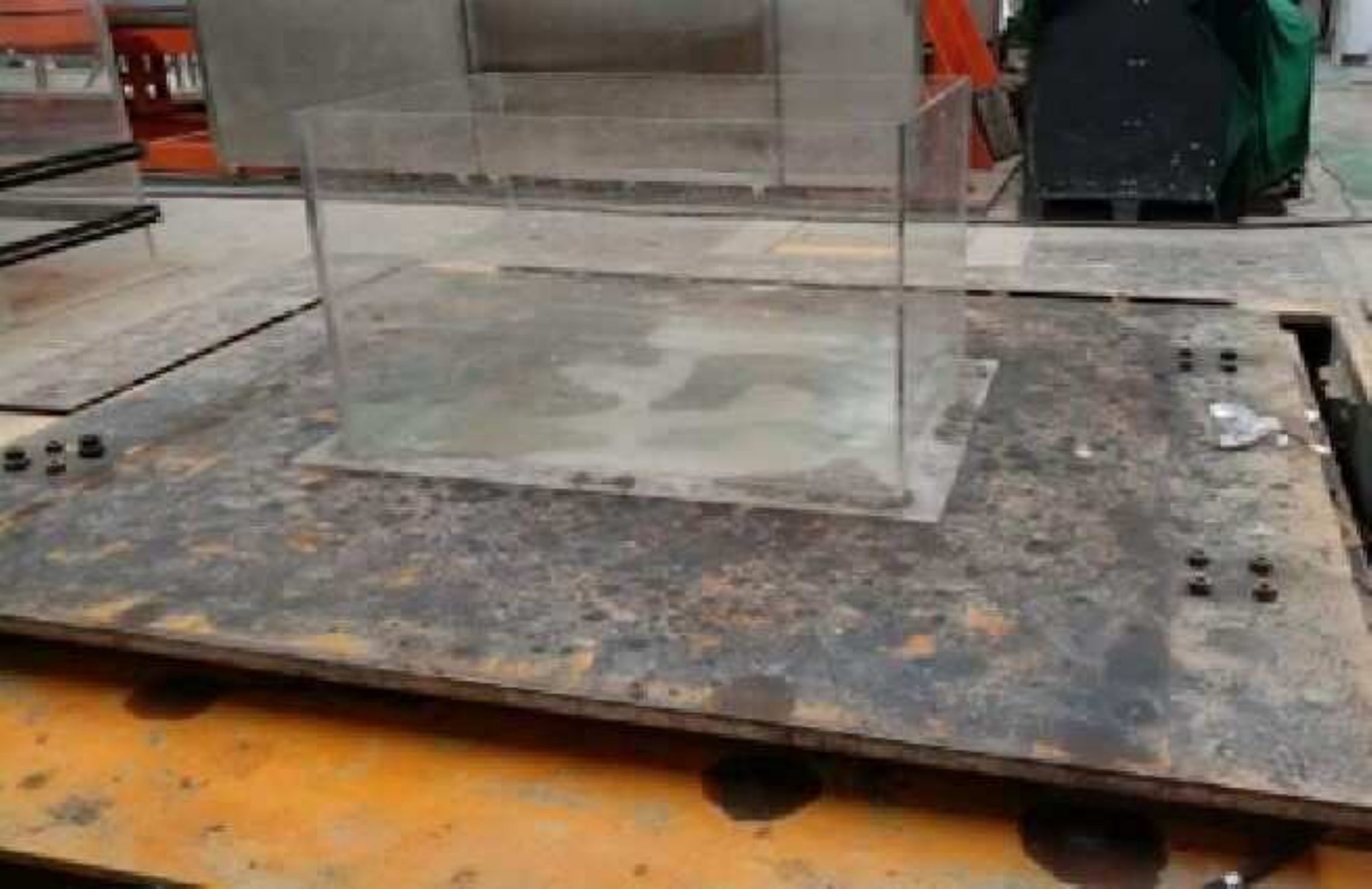}
          \end{center}
          \caption{\footnotesize{Laboratory experiment setup composed by a water tank installed above a shaking table.}}
          \label{fig:lab_experiment}
    \end{figure}

     Firstly, we adopt the real-time hybrid simulation technology \cite{tang2020performance}, which divides the test object into two parts: the physical substructure (i.e. TLD system) for which we would like to learn the model from the collected data; the numerical substructure which simulates the dynamics of a prescribed building. As shown in Fig. \ref{fig_tld}, after we input the earthquake signal under the structure, the relative displacement between the top of the structure and the TLD system can be transferred equivalently as the input of the shaking table so that we can measure the relative acceleration {\mz $u$} for the bottom of the water tank via the accelerometer. Then the reactive force $F$ of the water is considered as the output which is measured by the three-component force sensor.
     The setting of this test object  can be seen in Fig. \ref{fig:lab_experiment}: it is composed by a water tank of dimensions $0.8 \times 0.8 \times 0.268 $ $m^3$ (length $\times$ width $\times$ depth) placed above the shaking table.
    Moreover, the earthquake signals considered here are  two classical  seismic waves, namely the \textit{El Centro wave} and the \textit{Hachinohe wave}. Finally, the relative acceleration $u$ and reactive force $F$ signals gathered from this procedure were filtered in order to reduce the noise and then taken into account to build the datasets for learning the model {\mzzz parameters}. These two {\mzzz signals}, sampled with {\mzzz sampling} time $T_s = 0.001 s$, correspond respectively to the variables $u$ and $y$ introduced in Section \ref{sec: TDL_learning}.

    The evaluation of the Algorithm \ref{alg: robust EKF} was performed in MATLAB. In order to evaluate the effectiveness of the proposed approach, the simulations were done using both the EKF and the {\mzzz REKF}. The fixed model parameters are  $m_t = 171.520$ $kg$ and $\xi = 0.005$, while the reference values of the objects of interest drawn from the Housner model  are  $\beta = 0.612$ and $w = 5.489$ $rad / s$.
    Both datasets were tested in two different scenarios. The first one considers the initial conditions (I.C.) of  $\beta$ and  $\omega$ approximately known. The second case supposes instead that we do not have {\mzzz this a priori} knowledge, so they are set equal to the lower bounds. Accordingly, the I.C. of the predicted state vector for these frameworks are respectively:
\begin{align}
        \label{x01}x_{0} &= [0.01, -0.01, 0.5, 5]\\
        \label{x02} x_{0} &= [0.01, -0.01, 0.1, 1].
    \end{align}
    The output noise variance $R$ was assumed to be equal to $1$ for all the simulations. Regarding the value of the tolerance $c_t$, it has been designed as an exponential with decay rate of $0.01$ starting from the initial value $c_0 = 0.001$: $c_t=c_0e^{-0.001t}$. In plain word, we need robustness only at the early stage where the linearisation may be inaccurate due by the poor knowledge of the actual values for $x_3$ and $x_4$.
    In what follows, the initial state error estimation covariance matrix $V_0$, the process covariance matrix $Q = B B^T$ and the estimation results are going to be presented for the two datasets separately.

\begin{figure}[h]
            \centering
            \includegraphics[scale = 0.55,width=8.3cm,height=5.5cm]{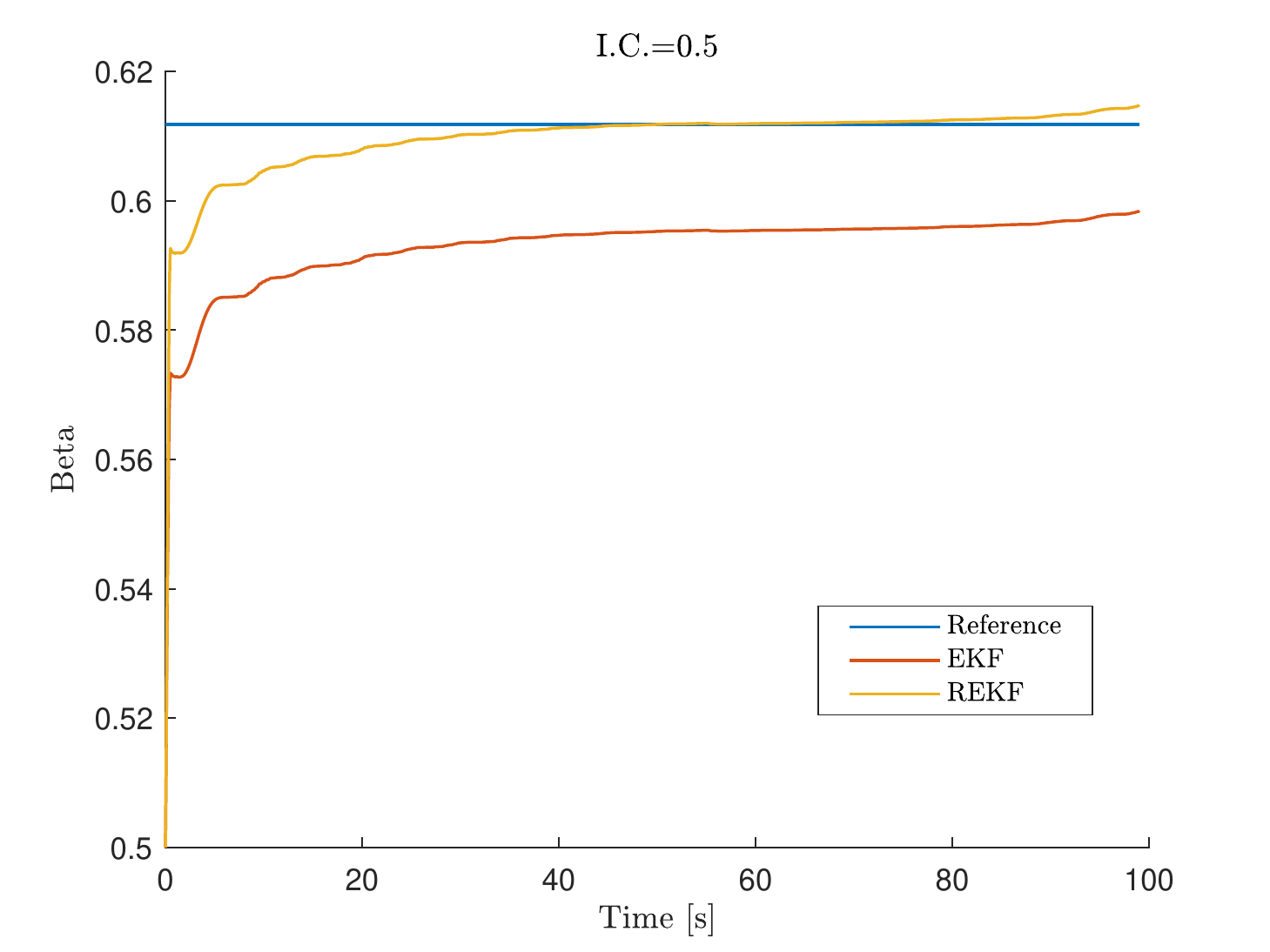}
            \includegraphics[scale = 0.55,width=8.3cm,height=5.5cm]{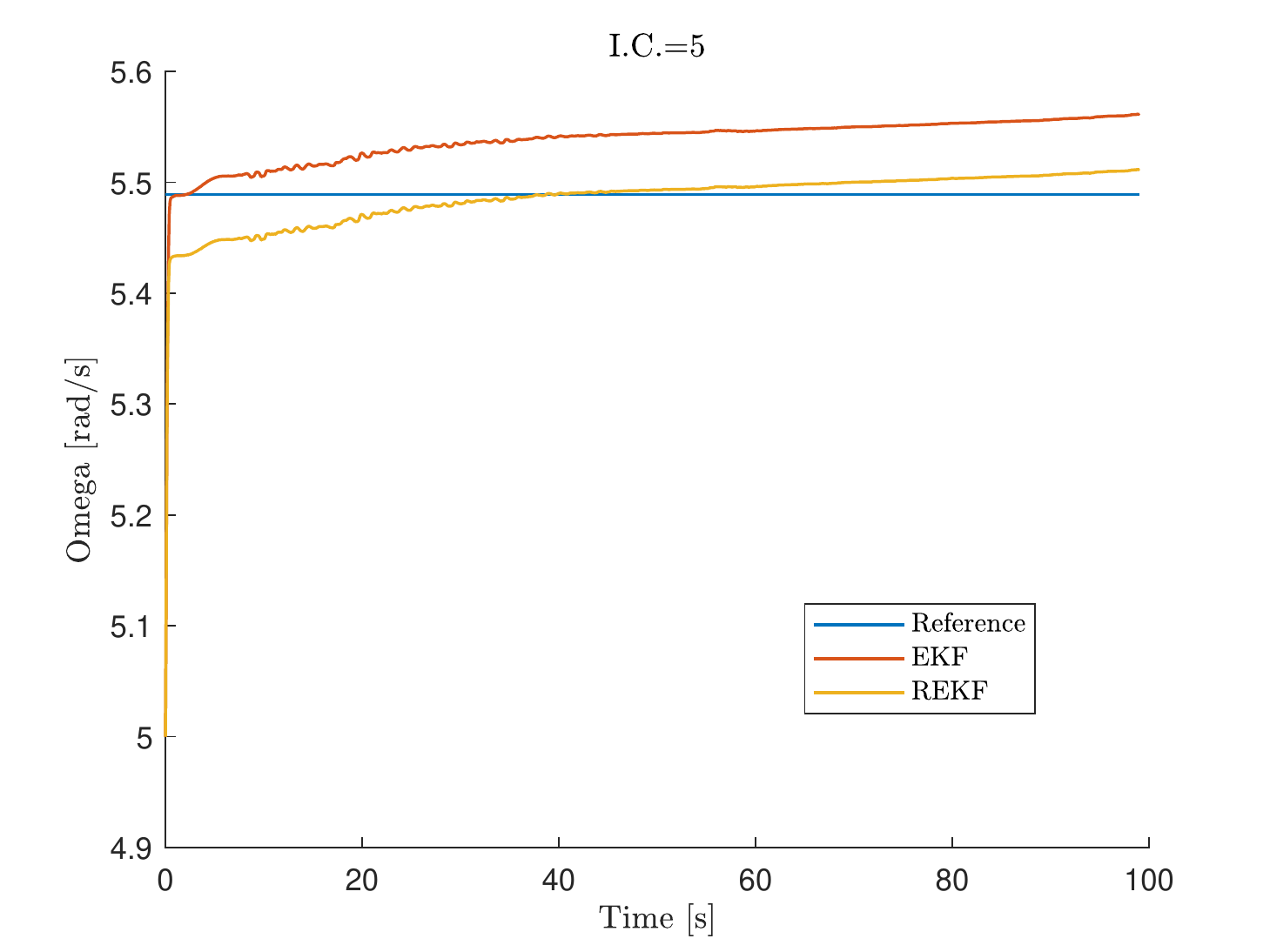}
            \caption{Estimation profiles of $\beta$ and $\omega$ given the initial conditions close to the reference ones (i.e. corresponding to (\ref{x01}) and (\ref{Q1})) and input $u$ measured from the \textit{El Centro wave}.}
            \label{fig: ElCentro_f1}
        \end{figure}

\begin{figure}[h]
            \centering
            \includegraphics[scale = 0.55,width=8.3cm,height=5.5cm]{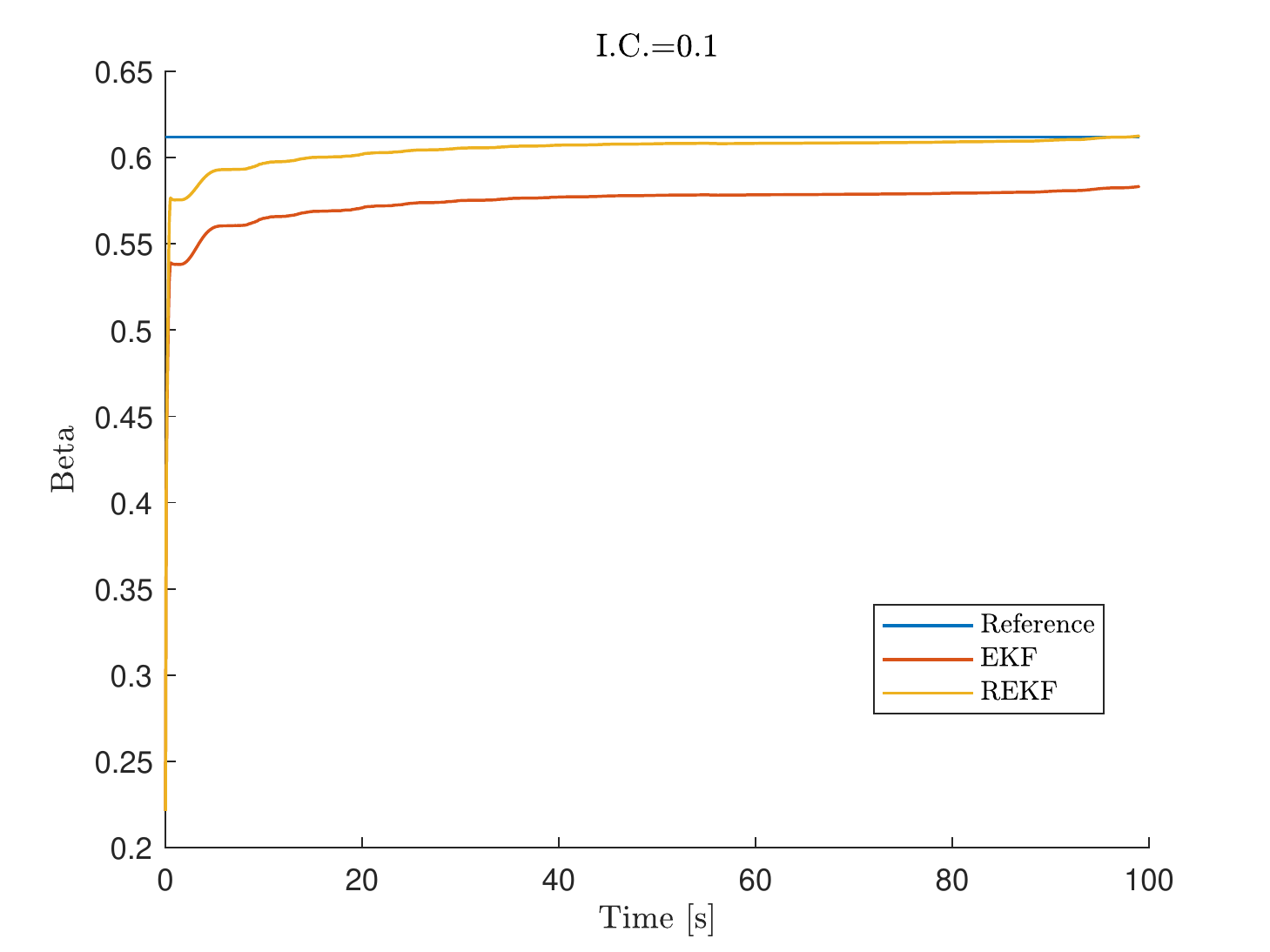}
            \includegraphics[scale = 0.55,width=8.3cm,height=5.5cm]{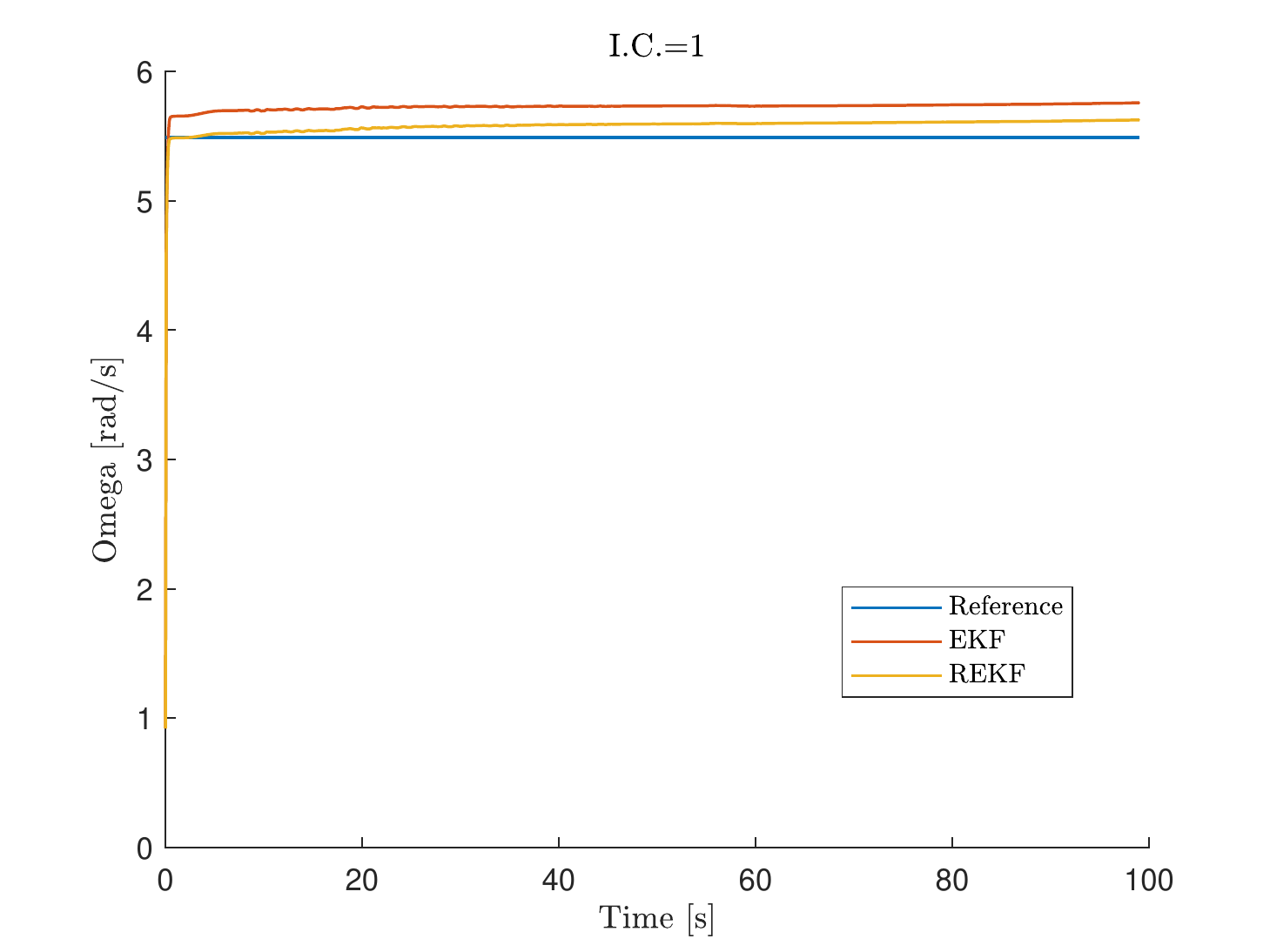}
            \caption{Estimation profiles of $\beta$ and $\omega$ given the initial conditions equal to their lower bounds (i.e. corresponding to (\ref{x01}) and (\ref{Q11})) and input $u$ measured from the \textit{El Centro wave}.}
            \label{fig: ElCentro_f2}
        \end{figure}

\textbf{El Centro wave dataset.}  The covariance matrix for the initial state is:
        \begin{equation*}
            V_0 = diag(1, 1, 0.001, 0.1).
        \end{equation*}
        For what concerns the process variance matrix, it was necessary to tune the variances of $\beta$ and $\omega$ depending on whether the values of the initial conditions of their estimated states were assumed to be close to the reference ones or not. The values of $Q$ used in the simulations for (\ref{x01}) and (\ref{x02}) were, respectively:        \begin{align}
\label{Q1}            Q &= diag(1, 1, 0.00115, 0.065),\\
 \label{Q2}           Q &= diag(1, 1, 0.0018, 0.145).
        \end{align}The results obtained from these initial states of the {\mzzz REKF} parameters can be seen respectively in Fig. \ref{fig: ElCentro_f1} and \ref{fig: ElCentro_f2}. It is worth noting that $\sigma_3$ and $\sigma_4$ in (\ref{Q1}),  corresponding to the I.C. for $\beta$ and $\omega$  close from the reference values, are smaller than the ones in (\ref{Q2}), corresponding to the I.C. for $\beta$ and $\omega$ far away from the reference values. {\sy In a plain word, these two choices hinge on this principle}: if the initial estimates of $\beta$ and $\omega$ are not close to the actual one, then we expect that the changing rate of the estimates of $\beta$ and $\omega$ is high. Such an a priori information is embedded in the stochastic hypermodel by choosing $\sigma_3$ and $\sigma_4$ ``sufficiently  large''.

\begin{figure}[h]
            \centering
            \includegraphics[scale = 0.55,width=8.3cm,height=5.5cm]{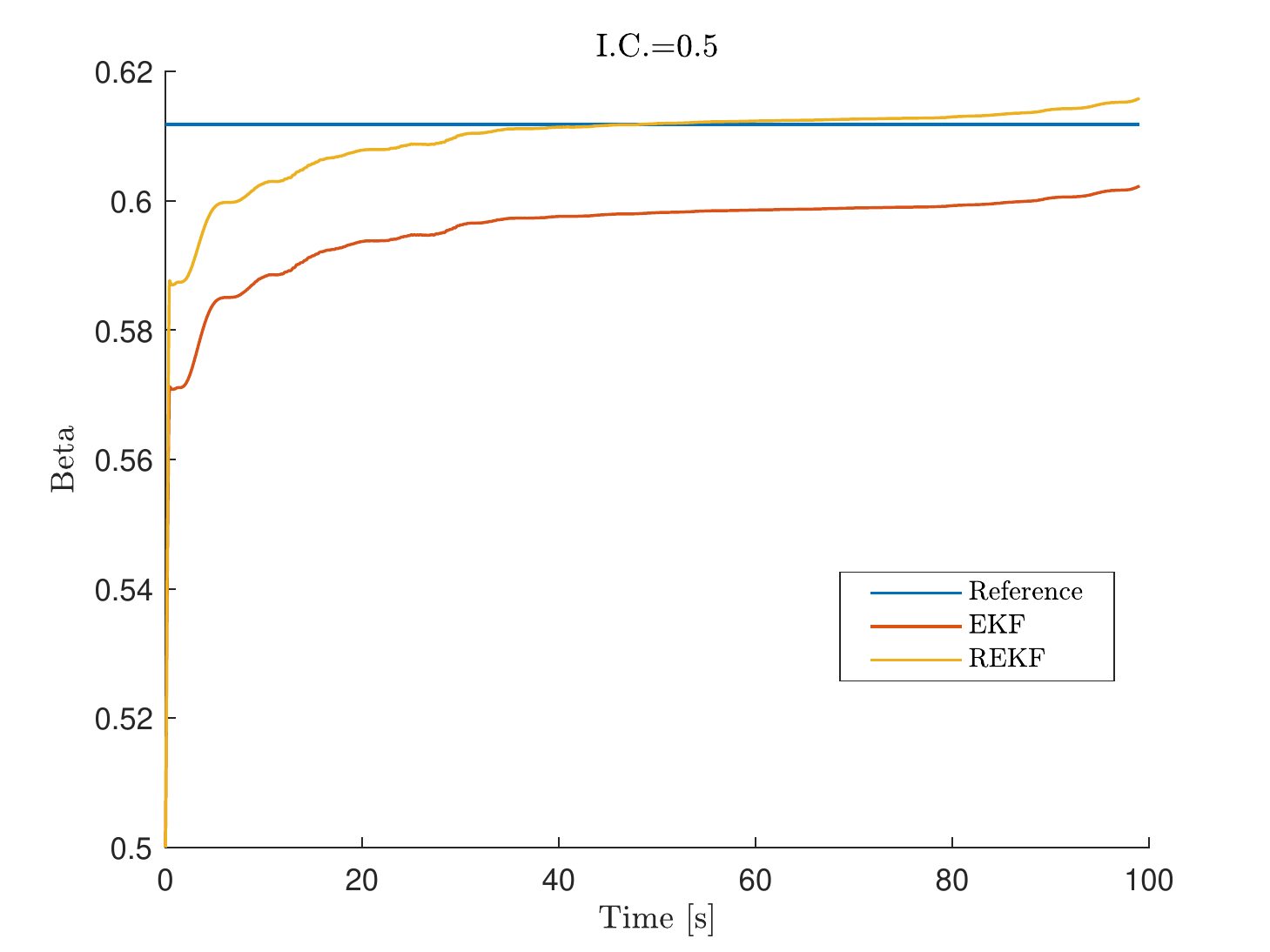}
            \includegraphics[scale = 0.55,width=8.3cm,height=5.5cm]{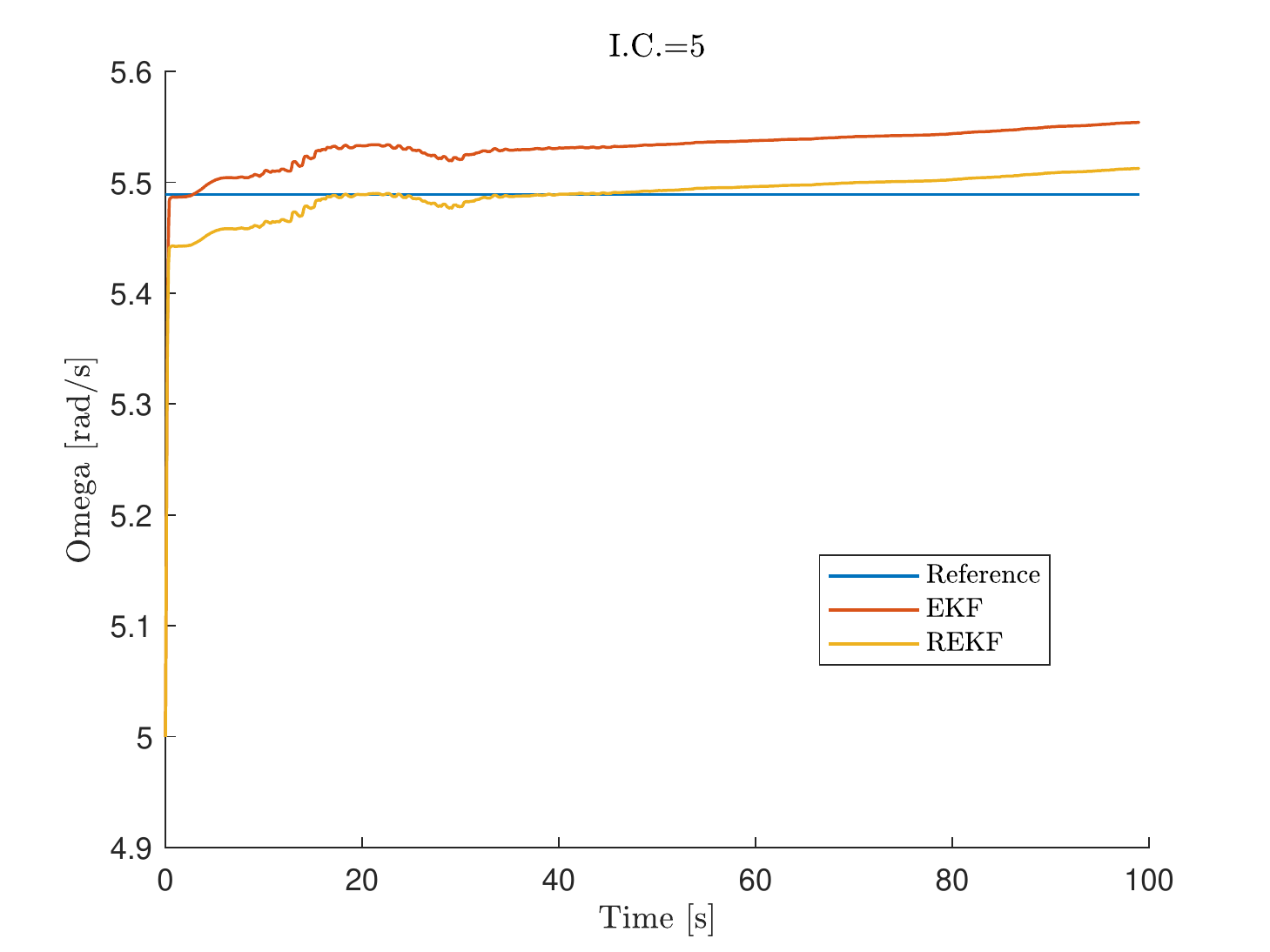}
            \caption{Estimation profiles of $\beta$ and $\omega$ given the initial conditions close to the reference ones (i.e. corresponding to (\ref{x02}) and (\ref{Q22})) and input $u$ measured from the \textit{Hachinohe wave}.}
            \label{fig: Hachinohe_f1}
        \end{figure}

\begin{figure}[h]
            \centering
            \includegraphics[scale = 0.55,width=8.3cm,height=5.5cm]{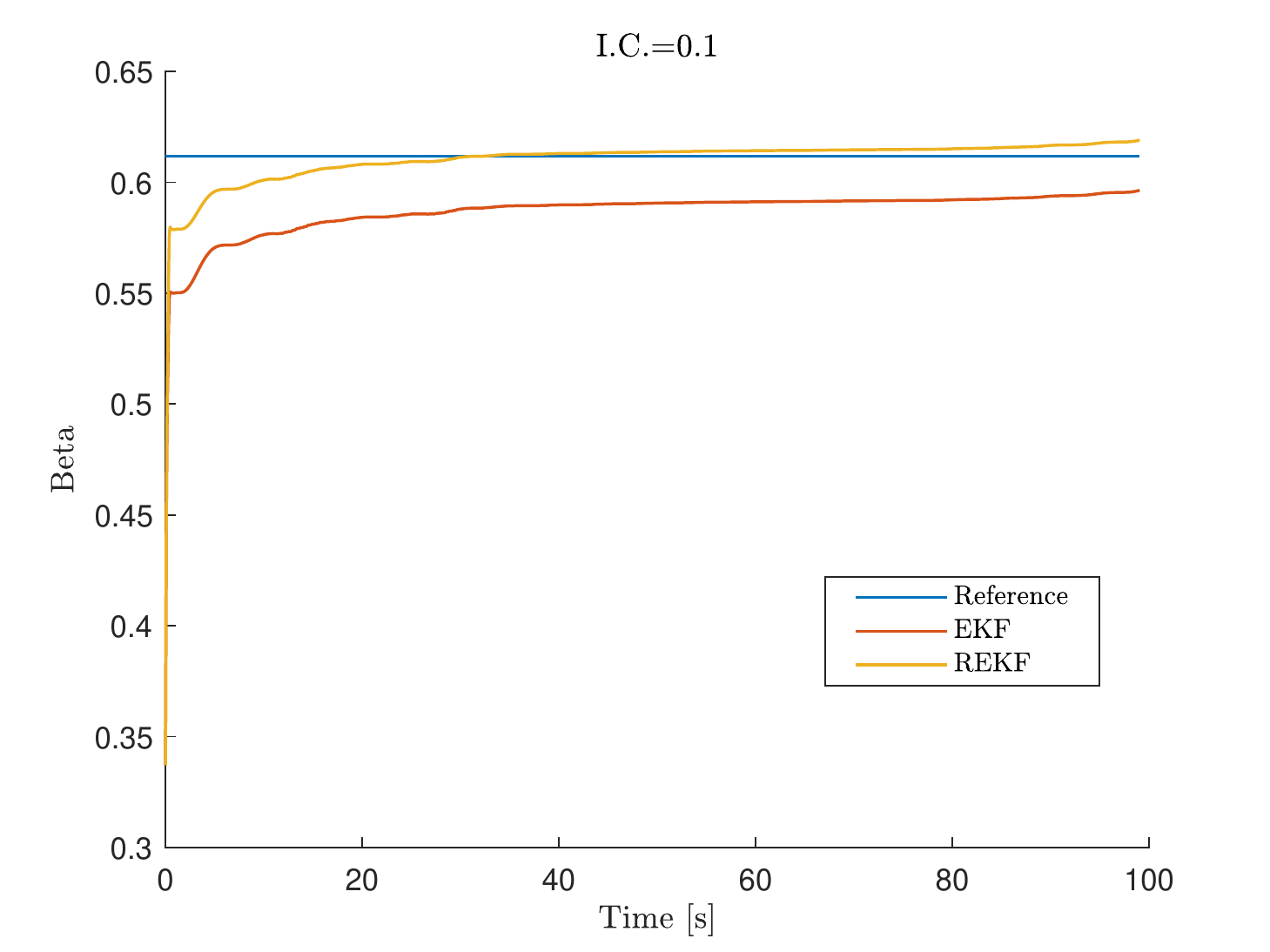}
            \includegraphics[scale = 0.55,width=8.3cm,height=5.5cm]{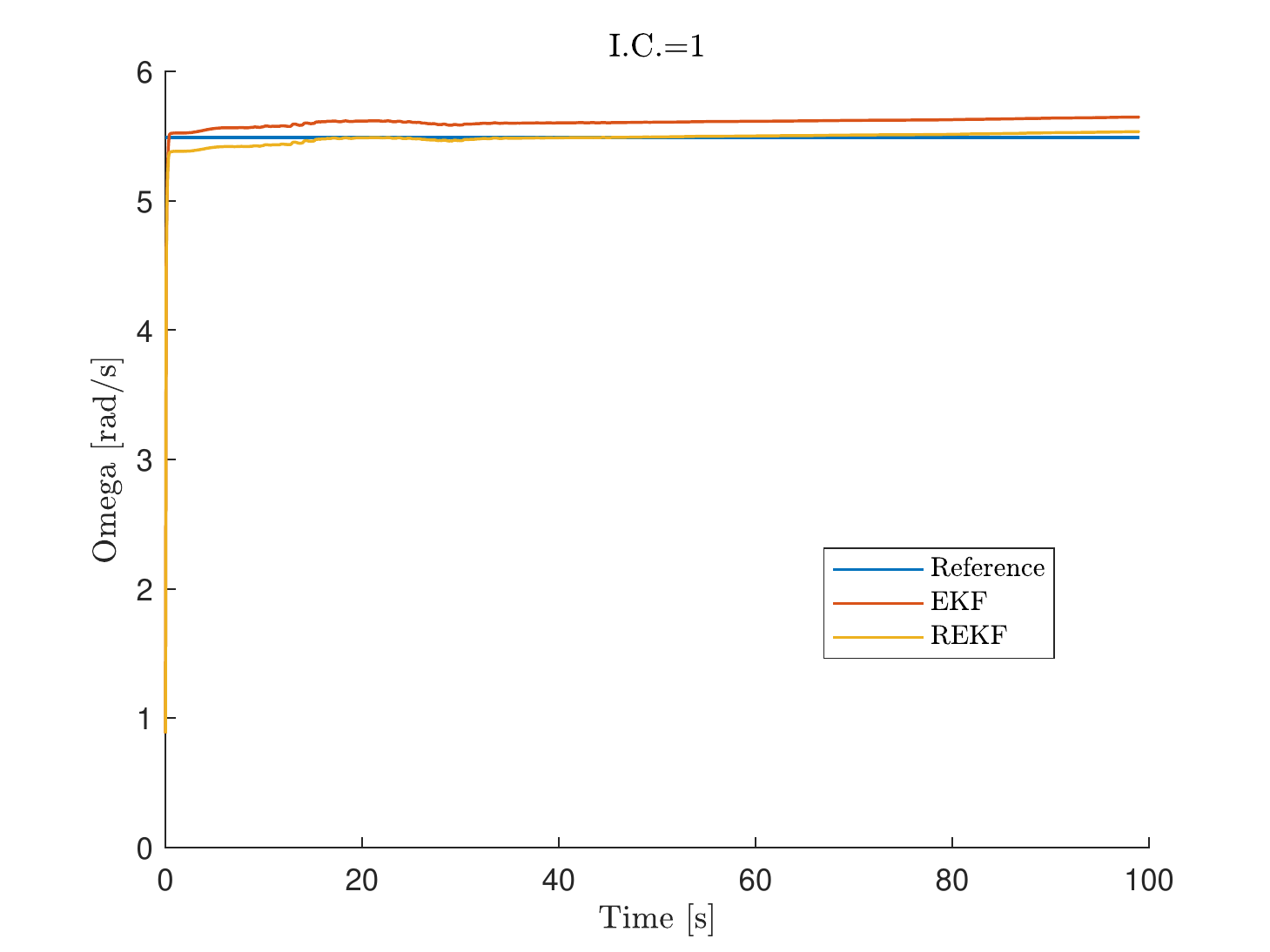}
            \caption{Estimation profiles of $\beta$ and $\omega$ given the initial conditions equal to their lower bounds (i.e. corresponding to (\ref{x02}) and (\ref{Q2})) and input $u$ measured from the \textit{Hachinohe wave}.}
            \label{fig: Hachinohe_f2}
        \end{figure}

    \textbf{Hachinohe wave dataset.} The covariance matrix of the initial state is:
        \begin{equation*}
            V_0 = diag(1, 1, 0.002, 0.15).
        \end{equation*}
        Following the same principle of before, we consider the following covariance matrices for (\ref{x01}) and (\ref{x02}), respectively:

        \begin{align}
     \label{Q11}       Q&= diag(1, 1, 0.0013, 0.07),\\
\label{Q22}
            Q&= diag(1, 1, 0.0018, 0.15).
        \end{align} Fig. \ref{fig: Hachinohe_f1} and \ref{fig: Hachinohe_f2} show  the estimation results obtained with the discussed initial conditions.

    From Fig. \ref{fig: ElCentro_f1}-\ref{fig: Hachinohe_f2} it is possible to observe how the {\mzzz REKF} is able to identify correctly the parameters $\beta$ and $\omega$ for both datasets. Moreover, the robust approach performs always better with respect to the standard one, as expected. Taking into account only the profiles of the {\mzzz REKF}, they always reach the reference values in at most $40 s$ out of $100 s$. Furthermore, once the true value has been reached, the estimate tends not to deviate from the reference parameter keeping the relative error not larger than $1\%$ for $\beta$ and $0.5\%$ for $\omega$. These relative errors for the standard EKF are $2\%$ for $\beta$ and $1\%$ for $\omega$. We conclude that, for both the initial conditions, the performance of the {\mzzz REKF} is satisfactory and outperforms the standard EKF.


 \section{Conclusion} \label{sec_5}
In this paper, we learn the TLD dynamics by means of a robust extended Kalman filter whose uncertainty is expressed incrementally. The aim of this robust approach is to alleviate the model uncertainty given by the discretization of the Housner model as well as the poor knowledge of the noise process of the state space model. The robust algorithm has been tested with experimental data whose earthquake signal takes classical seismic waves (i.e. the El Centro wave and the Hachinohe wave). {\mzzz It is worth noting that the computational time required by (the present implementation of) the proposed  approach is larger than that of the standard one. However,} the experimental results showed that  the proposed method not only achieve a satisfactory accuracy, but it is also able to deal with two different initial conditions. Although this is a promising result and since we have collected hundreds of datasets, we envision to perform a large number of experiments to verify the reliability of the proposed algorithm.


\section{Acknowledgments}
This work is partially supported by NSFC under Grant Nos. 51978016. The authors would like to thank Prof. Zhenyun Tang at Beijing University of Technology for his support in the experimental platform and data collection.

\end{document}